\documentclass{amsart}
\usepackage{amsmath,amsthm,amssymb,xypic}
\usepackage[all]{xy}

\title{Equivalent birational embeddings}

\author{Massimiliano Mella
\and Elena Polastri}
\address{Dipartimento di Matematica\\
Universit\`a di Ferrara\\
Via Machiavelli 35\\
44100 Ferrara, Italia} \email{mll@unife.it \ plslne@unife.it}

\date{July 2008}
\thanks{Partially supported by Progetto PRIN 2006 ``Geometria
  sulle variet\`a algebriche'' MUR}

\theoremstyle{plain}
\newtheorem{Th}{Theorem}
\newtheorem{Co}{Corollary}
\newtheorem{theorem}{Theorem}[section]

\newtheorem{lemma}[theorem]{Lemma}

\theoremstyle{definition}

\newtheorem{definition}[theorem]{Definition}

\theoremstyle{remark}

\newtheorem{remark}[theorem]{Remark}

\DeclareMathOperator{\Bsl}{Bs}

\DeclareMathOperator{\mult}{mult}

\DeclareMathOperator{\Bir}{Bir}

\newcommand{\QED}{\ifhmode\unskip\nobreak\fi\quad {\rm Q.E.D.}} 

\newcommand{\f}{\varphi}
\newcommand{\A}{\mathcal{A}}
\newcommand{\C}{\mathbb{C}}

\newcommand{\G}{\mathcal{G}}

\renewcommand{\H}{\mathcal{H}}

\newcommand{\I}{\mathcal{I}}

\renewcommand{\L}{\mathcal{L}}
\newcommand{\M}{\mathcal{M}}

\renewcommand{\O}{\mathcal{O}}
\renewcommand{\P}{\mathbb{P}}

\newcommand{\rat}{\dasharrow}

\newcommand\lag{\langle}
\newcommand\rag{\rangle}

\begin{document}

\maketitle

\begin{abstract} Let $X$ be a projective
  variety of dimension $r$ over an algebraically closed
  field. It is proven that two birational
  embeddings of $X$ in $\P^n$, with
  $n\geq r+2$ are
  equivalent up to Cremona
  transformations of $\P^n$.
\end{abstract}

\section*{Introduction}

Let us consider a rational variety $X\subset
\P^n$. Then there exists a birational map $\phi:X\rat\P^r$. The
simplest embedding of $\P^r$, as a projective variety, is the
linear one. It is quite natural to ask whether the map $\phi$ can
be extended to a birational map $\Phi:\P^n\rat \P^n$ such that
$\Phi(X)$ is linear. We shall say that in this case the variety $X$
is {\sl Cremona equivalent}, see Definition
\ref{def:cremona}, to a linear space. This extension property reminds us
the Abhyankar-Moh property (AMP), \cite{AM}. The latter asks about
extensions of polynomial embeddings in $\C^n$ to automorphisms of
$\C^n$. The AMP has been studied extensively and seems granted for
high codimension smooth varieties,
see for instance
\cite{Je}.  The AMP can be extended
to affine varieties over an infinite field,
say $k$. Then one can ask whether two
different embeddings of the same affine variety are
equivalent up to an automorphism of
$k^n$. Also in this context the
answer is  positive for varieties with
isolated singularities and roughly local dimension half the
embedding dimension, see \cite{Sr} for the
precise statement.

The aim of this paper is to solve a similar
question in the context of birational
geometry of projective varieties.
We want to understand when two birational
embeddings of the same variety are equivalent
up to a Cremona transformation of the
projective space, in this case we say that
they are Cremona equivalent, see Definition
\ref{def:cremona} for the precise statement. 
In the case of complex field this has been
studied in \cite[Theorems 2.1, 2.2]{Je} where
a positive answer 
is given for varieties of codimension greater than
the dimension.

The birational nature of the problem suggests
that singularities should not be a main
issue. Indeed we are able to treat arbitrarily
irreducible and
reduced projective varieties.
The main novelty comes from the range where
the answer is positive.
 It
is not difficult to give examples of rational
hypersurfaces that are not Cremona equivalent to a
hyperplane, see section \ref{sec:divisor} for
the details. In general birationally
equivalent divisors are not Cremona
equivalent. What is really surprising is that
this is the only case.

\begin{Th} Let $X$ be an irreducible and
  reduced projective variety of dimension $r$
  over an algebraically closed field $k$.
  Then two birational embeddings in $\P^n$ are
  Cremona equivalent as long as $n\geq r+2$.
\label{th:main}
\end{Th}

To prove the theorem we produce a chain of
Cremona transformations that modify the
linear systems giving the two
embeddings one into the other.
This is done by looking at the two
birational embeddings as different projections
of a common embedding.
We want to stress some by-products of the main
Theorem.
\begin{Co} Let $X\subset\P^n$ be a subvariety
  of codimension at least $2$.

If $X$ is rational then it is Cremona equivalent to a linear space. 

If $X$ is birational to a smooth
  subvariety of $\P^n$. Then it is possible to
  resolve the singularities of $X$ with a
  Cremona transformation.

The Cremona group of $\P^n$ contains
$\Bir(X)$, the group of birational
transformations of $X$.
\end{Co}

We like to look at the Theorem as a way
to say that the Cremona group of $\P^n$ is
really huge. We thank the referee for a
careful reading.


\section{Notations and preliminaries}
\label{sec:bir}

\noindent We work over an algebraically closed
field $k$.
\begin{definition} For any linear space
  $L\subset\P^n$ the map
$$\pi_L:\P^n\rat\P^{n-\dim L-1}$$ 
is the projection from $L$.
\end{definition}

Let $X$ be a projective irreducible and
reduced variety and $\L$ a
linear system on $X$. Assume that $\L$ is
generated by $\{L_0,\ldots,L_n\}$. Then the
map $\f_\L:X\rat \P(\L^*)$  is given by
evaluating the sections of $L_i$ at the points
of $X$. 
\begin{definition} We say that $(X,\L)$, or
  simply $\L$, is a
  birational embedding (in $\P^n$) if
  $\f_\L:X\rat \P^n$ is
  birational onto the image. We say that
  $\f_\L=\f_\M$, for two birational
  embeddings, if there exists a dense
open subset $U\subset X$ where $\f_\L$ and
$\f_\M$ are both defined and have equal restriction.
\end{definition}

\begin{remark}
\label{rem:div} Note that given a birational embedding
  $\L$ in $\P^n$ we can consider it also an
  embedding into $\P^{n+h}$ by adding
  $h$-times the zero section to get 
  $$\L'=\{L_0,\ldots, L_n,0,\ldots,0\}$$
In all that follows we apply, mainly without mention,
this construction  to
compare birational embeddings into different
projective spaces.

  We are
  interested in studying  birational embeddings
  of a fixed variety $X$. We therefore
  identify $\L$  with
   $\O_{\f_\L(X)}(1)$ via
  $(\f_\L)^{-1}_*(\O_{\P^n}(1))$.

 Let $D\subset X$ be a divisor and
consider the linear system
$\L_D=\{DL_0,\ldots,DL_n\}$. Then we have
$\f_\L=\f_{\L_D}$. In what follows we identify
$\L$ and $\L_D$.
\end{remark}

\begin{definition}
 \label{def:cremona}
Let  $\L$ in $\P^n$ and $\G=\{G_0,\ldots,G_r\}$ in $\P^r$ be two birational
  embeddings, assume that $n\geq r$ and let
  $\G'=\{G_0,\ldots,G_r,0,\ldots,0\}$,
  obtained by adding $(n-r)$-times the
  zero section. We say that $\L$ is Cremona
  equivalent to $\G$, or simply equivalent, if there
  exists a birational map $\Phi:\P^n\rat\P^n$
  such that
$$\f_\L=\Phi\circ\f_{\G'}$$
Such a $\Phi$ is called a (Cremona) equivalence between
$\L$ and $\G$.
\end{definition}

\begin{remark}
The relation introduced is an equivalence
relation on the      birational linear systems
of a fixed variety. We never
ask  the linear system neither to be complete
nor to be minimally generated by the $L_i$'s.

 The equivalence $\Phi$ has to be defined on the
general point of $\f_{\G'}(X)$.
\end{remark}

Let $\L$ and $\G$ be two linear systems on
 $X$. 
Then we have the following commutative
 diagram

\[
\xymatrix{
 & &  \P^{N}\ar@{.>}[drr]^{\pi_G}\ar@{.>}[dll]_{\pi_L} &&&
                                                                \\
\P^n&&X\ar@{.>}[ll]^{\f_\L}\ar@{.>}[rr]_{\f_\G}\ar@{.>}[u]_{\f_{\L+\G}}&&\P^n
 }
\]

where $\L+\G=\{L_iG_j\}$, with $i,j=0,\ldots,n$, and 
$$\pi_L\circ\f_{\L+\G}=\f_{\L_{G_0}}\qquad
 \pi_G\circ\f_{\L+\G}=\f_{\G_{L_0}}$$
In other terms  $L$ and $G$ are
linear spaces spanned, respectively by
 $\{G_jL_i\}$,  with
$j=1,\ldots, n$, $i=0,\ldots,n$,  and $\{G_jL_i\}$ with
$i=1,\ldots, n$, $j=0,\ldots,n$. 

This diagram allows us to look at two
different embeddings as projections
from a common one.

A natural way to construct a birational self-map
of $\P^n$ is to consider a hypersurface
of degree $k$
in $\P^{n+1}$ with two points of multiplicity
exactly $k-1$. Then the projections from the singular
points build up the required self-map. The
following  Lemma allow us to use this trick in
a wide contest.

\begin{lemma}
\label{le:hyper}
Let $Y\subset\P^{n+1}$ be an irreducible
reduced variety and $q_1$, $q_2$ two points in
$\P^{n+1}$.  Let $CY_i$ be
  the cone over $Y$ with vertex $q_i$.
Assume that $\dim
Y\leq n-2$,
$h^0(\I_{Y}(1))\neq 0$ and $CY_i\not\subset\Bsl| \I_Y(1)|$.
Then for $k\gg 0$ there exists an irreducible
reduced hypersurface $S\in|\I_Y(k)|$ with:
\begin{itemize}
\item[-]
$\mult_{q_i}S=k-1$, for $i=1,2$,
\item[-] $S\not \supset CY_i$.
\end{itemize}
\end{lemma}
\begin{proof}
Let $l=\lag q_1,q_2\rag$ be the
line spanned by the $q_i$'s and consider the projections
$$\pi_{q_1}:\P^{n+1}\rat\P^{n},\ \pi_{q_2}:\P^{n+1}\rat\P^{n},\ \pi_l:\P^{n+1}\rat\P^{n-1}$$
Let $\tilde{Y}=\pi_l(Y)$, $Y_i=\pi_{q_i}(Y)$
be varieties. Then we have the following
diagram
\[
\xymatrix{
&Y\subset\P^{n+1}\ar@{.>}[dl]_{\pi_{q_1}}\ar@{.>}[dr]^{\pi_{q_2}}\ar@{.>}[dd]_{\pi_{l}} &&&
                                                                  \\
Y_1\subset\P^n\ar@{.>}[dr]&&\P^n\supset Y_2\ar@{.>}[dl]\\
&\tilde{Y}\subset\P^{n-1}&
 }
\]

Let us consider $D=(d=0)\subset\P^{n-1}$ a hypersurface of
degree $\delta$ with 
$\tilde{Y}\subset D$. Let $H=(h=0)\in
|\I_{Y}(1)|$ be a general hyperplane.
Define
$$S=(dg_1+hg_2=0)\subset\P^{n+1} $$
where:
\begin{itemize}
\item[$g_1$]  is general of degree $k-\delta$ and
multiplicity $k-\delta-1$ at $p_i$;
\item[$g_2$]  is general with $\mult_{q_i}hg_2=k-1$.
\end{itemize}
It is easy to check that $S$ satisfies all the
requirements. 
\end{proof}

\section{The proof}
We are now ready to prove Theorem
\ref{th:main}.
Let $X$ be an irreducible reduced projective
variety. Let $\L$ and $\G$ be two birational
embeddings in $\P^n$.  Keep in mind
that they are both projections of
$\L+\G$ and
$$\L=\L_{G_0}=\{L_0G_0,\ldots,L_nG_0\},\ \ \ \  \G=\G_{L_0}=\{L_0G_0,\ldots,L_0G_n\}.$$
We want to
construct a sequence  of Cremona
equivalent 
linear systems $\{\A_i\}$, for
$i=0,\ldots n$, with
\begin{itemize}
\item[-] $\A_0=\L_{G_0}=\{L_0G_0,\ldots,L_nG_0\}$;
\item[-] $\A_i=\{L_0G_0,\ldots,L_0G_i, 
  A^i_{i+1},\ldots, A^i_n\}$,
for some
  $A^i_j\in|\L+\G|$, $j=\hbox{i+1},\ldots, n$;
\item[-]  $\A_n=\G_{L_0}=\{L_0G_0,\ldots,L_0G_n\}$.
\end{itemize}

To prove the theorem we give a recipe that builds
$\A_{i+1}$ from $\A_i$.
 Let $\H_i=\{\A_i,L_0G_{i+1}\}$ be a
  linear system and
  $\f_{\H_i}:X\to\P(\H_i^*)=\P^{n+1}$  the
  associated embedding. Then we have
  $X_i=\f_{\A_i}(X)\subset
  (x_{n+1}=0)\subset\P^{n+1}$. Let
  $Z\subset\P^{n+1}$ be
  the cone over $X_i$ with vertex $q_{1}=[0,\ldots,0,1]$
  and $Y_i:=Z\cap H$ a general hyperplane section of
  $Z$. Then we have
$$\pi_{q_{1}}(Y_i)=X_i$$
and $Y_i$ birational to $X$.
Let $q_2\in
(x_0=\ldots=x_i= x_{n+1}=0)$ be a general
point. In particular $q_2\not\in H$ hence the projection
$\pi_{q_2|Y_i}$ is birational onto the
image. Let $X_{i+1}:=\pi_{q_2}(Y_i)$, then
$X_{i+1}$ is birational to $X$  and we define
$\A_{i+1}:=\O_{X_{i+1}}(1)$, keep in mind
Remark \ref{rem:div}. The choice of
$q_2$  gives 
$$\A_{i+1}=\{L_0G_0,\ldots,L_{0}G_{i+1},A^{i+1}_{i+2},\ldots,A^{i+1}_n\}$$
for some elements $A^{i+1}_j\in|\L+\G|$.

This reads in the following diagram
\[
\xymatrix{
&Y_i\subset\P^{n+1}\ar@{.>}[dl]_{\pi_{q_{1}}}
 \ar@{.>}[dr]^ {\pi_{q_2}}&\\
 \f_{\A_{i}}(X)\subset\P^n&&\P^n\supset\f_{\A_{i+1}}(X)\\
}
\]
By construction $h^0(\I_{Y_i}(1))\neq 0$ and
$\Bsl|\I_{Y_i}(1)|\not\ni q_i$, for
$i=1,2$. By hypothesis we have $\dim Y\leq n-2$,
then by Lemma \ref{le:hyper} there exists an
irreducible hypersurface
$S\in|\I_{Y_i}(k)|$ with $\mult_{q_i}
S=k-1$ and not containing the cones over $Y_i$
with vertex both $q_1$ and $q_2$.
In particular
$$\pi_{q_{1}|S}\ {\rm and}\ \pi_{q_2|S}$$ 
are birational maps to $\P^n$ and
$\pi_{q_{1}|S}^{-1}$, $\pi_{q_{2}|S}^{-1}$ are defined on the general
point of $\f_{\A_i}(X)$,
$\f_{\A_{i+1}}(X)$ respectively.
Define the map $\Phi_i:\P^n\rat\P^n$ as follows
$$\Phi_i=\pi_{q_2|S}\circ \pi_{q_{1}|S}^{-1}$$
By construction $\Phi_i$ is an equivalence between $\A_i$ and $\A_{i+1}$.

\section{Divisorial embeddings}
\label{sec:divisor}
It is quite natural to expect that two
divisorial embeddings are in general not
Cremona equivalent. The argument of Noether--Fano
inequalities gives us a way to explicitly state it for
rational varieties, and
therefore show that Theorem \ref{th:main} is
the best possible result.
\begin{lemma} Let $X\subset\P^n$ be a rational
  variety of codimension $1$ and degree
  $d>1$. Assume that  $X$ is
  Cremona equivalent to a hyperplane. Then
  the singularities of the pair
  $(\P^n,(n+1)/d X)$ are not
  canonical.
\label{le:can}
\end{lemma} 
\begin{proof}
Let $\Phi:\P^n\rat\P^n$ be an equivalence and 
\[
\xymatrix{
&Z\ar[dl]_{p}
 \ar[dr]^q &\\
 \P^n\ar@{.>}[rr]^\Phi&&\P^n\\
}
\]
a resolution of $\Phi$.
 Then we have
$$\O_Z\sim
p^*(\O(K_{\P^n}+\frac{n+1}dX)=K_Z+\frac{n+1}d X_Z-\sum
a_iE_i$$
and
$$
q^*(\O(K_{\P^n}+\frac{n+1}d\O(1))=K_Z+\frac{n+1}d X_Z-\sum
b_iF_i $$
where $E_i$, respectively $F_i$, are $p$,
respectively  $q$, exceptional
divisors.
Let $l\subset\P^n$ be a general line in the
right hand side $\P^n$. Then we have
$$0 > q^{-1}l\cdot (q^*(\O(K_{\P^n}+\frac{n+1}d\O(1))+\sum
b_iF_i)=(\sum a_iE_i)\cdot q^{-1}l$$
This proves that at least one $a_i<0$ proving
the claim.
\end{proof}

Lemma \ref{le:can} allows to produce examples
of divisor that are not Cremona
equivalent to a hyperplane. 
Let $C\subset\P^2$ be a rational curve with
only ordinary double points.  If $\deg C\geq 6$
then $C$ is never Cremona equivalent to a
line. It is easy to produce such curves by
projecting a divisor of type $(1,a)$ in a
smooth quadric. 

One can construct similar examples in
arbitrary dimension.

\end{document}